\documentclass[12pt,notitlepage,twoside]{article}
\usepackage{latexsym}
\usepackage{amsmath}
\usepackage{amsfonts}
\usepackage{amssymb}
\renewcommand{\a }{\alpha }

\renewcommand{\d}{\delta }

\newcommand{\D }{\Delta }

\newcommand{\e }{\varepsilon }

\renewcommand{\l }{\lambda }
\renewcommand{\L }{\Lambda }

\newcommand{\n }{\nabla }
\newcommand{\var }{\varphi }

\newcommand{\Sig }{\Sigma}

\renewcommand{\o }{\omega }

\newcommand{\ov}{\overline}

\newcommand{\be}{\begin{equation}}
\newcommand{\ee}{\end{equation}}

\newenvironment{pfn}[1]{\noindent{\bf Proof of {#1}\enspace}}{
\hfill$\Box$\medskip}
\newcommand{\R}{\mathbb{R}}
 
\newcommand{\N}{\mathbb{N}}

\newtheorem{thm}{Theorem}[section]
\newtheorem{pro}[thm]{Proposition}
\newtheorem{lem}[thm]{Lemma}

\newtheorem{cor}[thm]{Corollary}

\numberwithin{equation}{section}

\author{Hichem CHTIOUI$^a$ \& Khalil EL MEHDI$^{b,c}$ \thanks{Corresponding author.
Present adress : see adress c mentioned above} \\
{\footnotesize
a : D{\'e}partement de Math{\'e}matiques, Facult{\'e} des Sciences de Sfax, Route
Soukra,}\\
{\footnotesize
 Sfax, Tunisia. E-mail : \texttt{Hichem.Chtioui@fss.rnu.tn} }\\
{\footnotesize
b :  Facult\'e des Sciences et Techniques, Universit\'e de Nouakchott, BP 5026, Nouakchott,}\\{\footnotesize
  Mauritania. E-mail : \texttt{khalil@univ-nkc.mr}}\\
{\footnotesize
 c : The Abdus Salam International Centre for Theoretical Physics, Mathematics Section,}\\
{\footnotesize
 Strada Costiera 11, 34014 Trieste, Italy. E-mail : \texttt{elmehdik@ictp.trieste.it}}
}

 \title { \Large \textbf{Prescribed Scalar Curvature  with Minimal \\
  Boundary Mean Curvature on $S^4_+$}}

\begin{document}

\date{ }

\maketitle

{\footnotesize
\noindent
{\bf Abstract.}
This paper is devoted to the  prescribed  scalar curvature under minimal boundary mean curvature on the standard four dimensional half sphere. Using topological methods from the theory of critical points at infinity, we prove some existence results.
}

\section{Introduction }
\mbox{}
In this paper, we are interested in some nonlinear equation arising from a geometric context. Namely, let
 $$
S_+^n =\{x=(x_1,...,x_{n+1})\in \R^{n+1}\, /\, |x|=1,\, x_{n+1}>0\}
$$
 be the standard half sphere endowed with its standard metric $g$, $n\geq 3$ and let \\ $f: S^n_+\to \R$ be a given function and we consider the following problem : does there exist a metric $\tilde{g}$ conformally equivalent to $g$ such that $R_{\tilde{g}} \equiv f$ and $h_{\tilde{g}} \equiv 0$ ?  where
$R_{\tilde{g}}$ is the scalar curvature of $S^n_+$ and $h_{\tilde{g}}$ is the
 mean curvature of $\partial S^n_+$, with respect to $\tilde{g}$. Setting $\tilde{g}=u^{4/(n-2)}g$  a conformal metric to $g$,
where $u$ is a smooth positive function, the above problem has the following analytical formulation : find a smooth positive function which solves the  following problem
\begin{eqnarray*}
(P)\quad \left\{
\begin{array} {ccccc}
L_gu := -\D _gu+ \frac{n(n-2)}{4}u & = & K u^{\frac{n+2}{n-2}} &\mbox{ in }& S_+^n\\
\frac{\partial u}{\partial \nu}&=&0&
\mbox{ on }& \partial S_+^n
\end{array}
\right.
\end{eqnarray*}
where $\nu$ is the outer normal vector with respect to $g$ and where $K=\frac{n-2}{4(n-1)}f$.\\
Such  kind of problem $(P)$ has attracted much attention
  (see \cite{ALM} \cite{BH}, \cite{BEO1}, \cite{BEO2}, \cite{C},\cite{DMO}  \cite{E1}, \cite{E2},  \cite{HL1}, \cite{HL2} \cite{Li} and the references therein).\\
The main difficulty one encounters in problem $(P)$ appears when we consider it from a variational viewpoint, indeed, the Euler functional associated to $(P)$ does not satisfy the Palais-Smale condition, that is, there exist noncompact sequences along which the functional is bounded and its gradient goes to zero. This fact is due to the presence of the critical Sobolev exponent in $(P)$. Moreover, there are topological obstructions to solving $(P)$, based on the Kazdan-Warner type condition, see \cite{bp}. Hence it is not expectable to solve problem $(P)$ for all functions $K$, and so a natural question arises :  under which conditions on $K$, $(P)$ has a positive solution?\\
 Yanyan Li \cite{Li} and Djadli, Malchiodi and Ould Ahmedou  \cite{DMO} studied problem $(P)$ on the three dimensional standard half
sphere. Their method involves a fine blow up analysis of some subcritical
approximations and the use of the topological degree tools.
 Ben Ayed, El Mehdi and Ould Ahmedou \cite{BEO1}, \cite{BEO2} gave some sufficient topological conditions on $K$ to find solutions to $(P)$ for $n$ bigger than or equal to $4$. Their approach uses algebraic topological tools from the theory of critical points at infinity (see Bahri \cite{B1}).\\
Notice that problem $(P)$ is, in some sense, related to the well-known Scalar Curvature problem on $S^n$
$$
(P') \qquad -\D_g u + \frac{n(n-2)}{4}u = K u^{(n+2)/(n-2)} \qquad \mbox{in } S^n$$
to which many works have been devoted (see for example the monographes \cite{A}, \cite{H} and  references therein.)\\
Regarding problem $(P')$, there is a difference between the three cases $n=3$, $n=4$ and higher dimensions. In the case $n=3$, the interaction between two of the functions failing the Palais-Smale condition dominates the self interaction, while for $n=4$, there is a balance phenomenon, and for $n\geq 5$ the self interaction dominates the interaction of two of those functions (see \cite{B2}, \cite{BaC}, \cite{BCCH}).\\
For problem $(P)$, such a balance phenomenon (i.e. the self interaction and the interaction are of the same size) appears in dimensions $3$ and $4$, see \cite{DMO}, \cite{BEO2}.
In this work, we focus on the four dimensional case to give more existence results. We are thus reduced to find positive solutions of the following problem
\begin{eqnarray*}
(1)\quad \left\{
\begin{array} {ccccc}
 -\D _gu+ 2u & = & K u^3 &\mbox{ in }& S_+^4\\
\frac{\partial u}{\partial \nu}&=&0&
\mbox{ on }& \partial S_+^4
\end{array}
\right.
\end{eqnarray*}
 Precisely, we borrow some of the ideas developed in  Bahri \cite{B2}, Aubin-Bahri \cite{AB}, Ben Ayed-Chtioui-Hammami \cite{BCH}, Ben Ayed-El Mehdi-Ould Ahmedou \cite{BEO1}, \cite{BEO2} and Chtioui \cite{HC}. The main idea is to precise the topological contribution of the critical points at infinity between the level sets of the associated Euler functional and the main issue is under our conditions on $K$, there remains some difference of topology which is not due to the critical points at infinity and therefore the existence of a solution of $(1)$.\\
In order to state our  results, we fix some notations and  assumptions that we are using in our results.

Let $G$ be the Green's function of $L_g$ on $S^4_+$ and $H$ its regular part defined by
\begin{eqnarray*}
\begin{cases}
G(x,y)=(1-\cos(d(x,y)))^{-1}-H(x,y),\\
\D H =0 \mbox{ in } S^4_+,\quad \partial G/\partial\nu =0 \mbox{ on } \partial S^4_+
\end{cases}
\end{eqnarray*}
 Let $K$ be a $C^3$ positive function on $\ov{S^4_+}$.\\

 Throughout this paper, we assume that  the following two assumptions hold\\
(1.1)\hskip 0.5cm  $K$ has only nondegenerate critical points
$y_0, y_1,...,y_m$ such that $y_0$ is the unique absolute maximum
of $K$ on $\ov{S^4_+}$ and such that
$$
\frac{-\D K(y_i)}{3K(y_i)} +4 H(y_i,y_i) \neq 0, \qquad \mbox{ for
} \, i=0,1,...,m
$$
(1.2)\hskip 0.5cm All the critical points of
$K_1=K_{/\partial S^4_+}$ are $z_1,...,z_{m'}$, and satisfy
$$
\frac{\partial K}{\partial \nu}(z_i) < 0, \qquad \mbox {for } i=1,...,m'
$$
Now we introduce the following set
$$
\mathcal{F}^+ =\{ y\in S^4_+/ \n K(y)=0 \,\, \mbox{ and
}\,\,\frac{-\D K(y_i)}{3K(y_i)} +4 H(y_i,y_i) > 0 \}
$$

Thus we  are able to state our first result
\begin{thm}\label{t:11}
If $y_0 \notin \mathcal{F}^+$, then
  problem $(1)$ has a solution.
\end{thm}
In the above result, we have assumed that $y_0 \notin
\mathcal{F}^+$. Next we want to give some existence result for problem $(1)$ when $y_0 \in \mathcal{F}^+.$ To this aim, we introduce some notation.\\
For $s\in\N^*$ and for any $s$-tuple $\tau_s=(i_1,...,i_s) \in (\mathcal{F}^+)^s$ such that $i_p\neq i_q$ if $p\neq q$, we define a Matrix $M(\tau_s)=(M_{pq})_{1\leq p,q\leq s}$, by
$$
M_{pp}=\frac{-\D K(y_{i_p})}{3K(y_{i_p})^2}+4\frac{
H(y_{i_p},y_{i_p})}{K(y_{i_p})}, \qquad M_{pq}
=-\frac{4G(y_{i_p},y_{i_q})}{\left(K(y_{i_p})K(y_{i_q})\right)^{1/2}}
\quad \mbox{for } p\neq q,
$$
and we denote by $\rho (\tau _s)$ the least eigenvalue of $M(\tau _s)$.
It was first pointed out by Bahri \cite{B1} (see also \cite{BaC}
and \cite{BCCH}), that when the self interaction and the
interaction between different bubbles are the same size,
the function $\rho$ plays a fundamental role in the existence
of solutions to problems like $(P)$. Regarding problem $(P)$,
Djadli-Malchiodi-Ould Ahmedou \cite{DMO} observed that such
kind of phenomenon appears when $n=3$.\\
Now let $Z$ be a pseudogradient of $K$ of Morse-Smale type (that
is the intersections of the stable and the unstable manifolds of
the critical points of $K$ are transverse).\\
(1.3) \hskip 0.5cm We assume throughout this paper that
 $W_s (y_i)\cap W_u(y_j) = \emptyset $ for any $y_i\in \mathcal{F}^+$
 and for any $y_j\notin \mathcal{F}^+ $,  where  $W_s (y_i)$ is the stable manifold of $y_i$ and $W_u(y_j)$ is the unstable manifold of $y_j$ for $Z$.\\
${\bf (H_1)}$\hskip 0.3cm Assume that $y_0 \in \mathcal{F}^+$\\
 Let $y_{i_1} \in \mathcal{F}^+ \diagdown \{y_0\}$ such that\\
(1.4) \hskip 2cm $K(y_{i_1}) = \max \{K(y)/y\in \mathcal{F}^+
\diagdown \{y_0\}\}$\\
 and we denote by $k_{i_1}=4-i(y_{i_1})$, where  $i(y_{i_1})$ is the Morse index of $K$ at  $y_{i_1}$.\\
 ${\bf (H_2)}$\hskip 0.3cm Assume that $i(y_{i_1}) \leq 3$.\\
We then have the following result :
\begin{thm}\label{t:12}
Under assumptions $(H_1)$ and $(H_2)$, if the following three
conditions hold
\begin{align*}
(A_0)& \quad  M(y_0,y_{i_1}) \,\,\mbox {is nondegenerate}\\
(A_1)& \quad \rho (y_0,y_{i_1}) <\quad 0\\
(A_2)& \quad \frac{1}{K(y)} > \frac{1}{K(y_0)} + \frac{1}{K(y_{i_1})} \qquad \forall y \in \mathcal{F}^+ \diagdown \{y_0, y_{i_1}\},
\end{align*}
then problem $(1)$ has a solution of Morse index $k_{i_1}$ or $k_{i_1}+1$.
\end{thm}
In contrast to Theorem \ref{t:12}, we have the following results based on a topological invariant for some Yamabe type problems introduced by Bahri \cite{B2}. To state these results, we need to fix  assumptions that we are using and some notation.\\
${\bf (H_3)}$\hskip 0.3cm Assume that $\rho (y_0,y_{i_1}) > 0$\\
Let
$$
X=\overline{W_s(y_{i_1})}
$$
Under (1.3) and (1.4), we derive that  $X= W_s(y_{i_1})\cup W_s(y_0)$. Thus
$X$ is a manifold of dimension $k_{i_1}$ without boundary.\\
We denote by $C_{y_0}(X)$ the following set
$$
C_{y_0}(X)=\{\a \d _{y_0}+(1-\a )\d _x \, / \, \a \in [0,1],\,
x\in X \},
$$
where $\d_x$ is the Dirac measure at $x$.\\
 For $\l $ large enough,
we introduce the map $f_\l : C_{y_0}(X)\to \Sigma ^+$, defined by
$$
 (\a \d _{y_0}+(1-\a )\d _x)  \longrightarrow
\frac{\a \d _{(y_0,\l )} +(1-\a )\d _{(x,\l )}}{||\a \d _{(y_0,\l
)}+(1-\a )\d _{(x,\l )}||},
$$
where $||.||$, $\Sig^+$ and $\d_{(x,\l)}$ are defined in the next section by \eqref{e:o1}, \eqref{e:o2} and \eqref{e:o3} respectively.\\
Notice that $C_{y_0}(X)$ and $f_\l ( C_{y_0}(X))$ are manifolds in dimension
$k_{i_1}+1$, that is, their singularities arise in dimension $k_{i_1}-1$ and
lower, see \cite{B2}. We observe that
 $ C_{y_0}(X)$ and $f_\l ( C_{y_0}(X))$ are contractible while $X$ is
not contractible.\\
For $\l $ large enough, we also define the intersection number(modulo 2) of
$f_\l (C_{y_0}(X))$ with $ W_s(y_0,y_{i_1})_\infty$
$$
\mu (y_0)=f_\l (C_{y_0}(X)). W_s(y_0,y_{i_1})_\infty ,
$$
where $W_s(y_0,y_{i_1})_\infty$ is the stable manifold of the
critical point at infinity $(y_0,y_{i_1})_\infty$ (see Corollary
\ref{c:32} below) for a decreasing pseudogradient $V$ for the Euler
functional associated to $(1)$ which is transverse to
$f_\l (C_{y_0}(X))$. Thus this number is well defined, see \cite{M}.\\
${\bf (H_4)}$\hskip 0.3cm Assume that $K(y_0) > 2K(y_{i_1})$. \\
We then have the following result
\begin{thm}\label{t:13}
Under  assumptions $(H_2)$, $(H_3)$ and $(H_4)$,
if $\mu (y_{i_1})=0$ then problem $(1)$ has a solution of Morse index $k_{i_1}$ or $k_{i_1}+1$.
\end{thm}
Now we give a statement  more general than Theorem \ref{t:13}.
To this aim, let $k\geq 1$, and define $X$ as the following
$$
X=\overline{\cup_{y\in B_k} W_s(y)}, \quad \mbox{ with } B_k \,\,\mbox{is any subset in }\,\,\{y \in \mathcal{F}^+ /
\, \mbox{i}(y)=4-k \},
$$
where i$(y)$ is the Morse index of $K$ at $y$.\\ 
${\bf (H_5)}$\hskip 0.3cm
We assume that $X$ is a stratified set without boundary (in the
topological sense, that is, $X\in \mathcal{S}_k(S^4_+)$,
the group of chains of dimension $k$ and $\partial X= 0$).\\
${\bf (H_6)}$\hskip 0.3cm Assume that for any critical point
$z$ of $K$ in $X\diagdown \{y_0\}$, we have $\rho (y_0,z) > 0$.\\
For $y\in B_k$ we define, for $\l $ large enough, the intersection
number(modulo 2)
$$
\mu (y)=f_\l (C_{y_0}(X)). W_s(y_0,y)_\infty
$$
By the above arguments, this number is well defined, see \cite{M}. \\
${\bf (H_7)}$\hskip 0.3cm Assume that $K(y_0) >2K(y)$
$\forall y\in \mathcal{F}^+ \diagdown \{y_0\}$.\\
Then we have the following theorem
\begin{thm}\label{t:14}
Under  assumptions $(H_5)$, $(H_6)$ and $(H_7)$,
if $\mu_(y)=0$ for any $y\in B_k$, then problem $(1)$ has a solution of Morse index $k$
or $k+1$.
\end{thm}
 Next we state a
 perturbative result for problem $(1)$. To this aim, we set
$$
 X=\overline{\cup_{y\in \mathcal{F}^+} W_s(y)}
$$
${\bf (H_{8})}$\hskip 0.3cm Assume that $X$ is not contractible
and denote by $m$ the dimension of the first nontrivial reduced homology group.\\
We then have
\begin{thm}\label{t:15}
Assume that 
   assumption $(H_8)$ holds. Thus there exists a constant $c_0$
independent of $K$ such that if 
$$ ||K-1||_{L^\infty (S^4_+)} \leq
c_0,$$
then problem $(1)$ has a solution of Morse index $\geq m$.
\end{thm}
Lastly, under the assumption  $(H_8)$, we can also find the following existence result :
\begin{thm}\label{t:16}
Under  assumption   $(H_8)$, if  the following two conditions hold
\begin{align*}
(C_1)& \quad  \mbox{for any }\,s,\,\,M(\tau_s) \,\,\mbox {is nondegenerate}\\
(C_2)& \quad \rho (y_i,y_j) < 0 \quad \forall y_i,\, y_j \in \mathcal{F}^+ \,\, \mbox{such that}\,\, y_i\neq y_j,
\end{align*}

then problem $(1)$ has a solution of Morse index $\geq m$.
\end{thm}
The rest of the paper is organized as follows.
In section 2, we set up the variational structure and recall
some known facts. Lastly,  section 3 is devoted to  the proofs
of our results.

\section { Some Known Facts }
\mbox{}
In this section we recall the functional setting and the
variational problem associated to $(1)$. We will also
recall some useful previous results.\\
Problem $(1)$ has a variational structure, the functional being
$$
J(u)= \frac {\int_{S_+^4}|\n u|^2+2\int_{S_+^4}u^2}
{\left(\int_{S_+^4}Ku^4\right)^{\frac{1}{2}}},
$$
defined on the unit sphere of $H^1(S^4_+)$ equipped with the norm
\begin{eqnarray}\label{e:o1}
||u||^2=\int_{S_+^4}|\n u|^2+ 2 \int_{S_+^4}u^2.
\end{eqnarray}
Problem $(1)$ is equivalent to finding the critical points of $J$ subjected to the constraint $u\in \Sig^+$, where
\begin{eqnarray}\label{e:o2}
\Sigma ^+=\{u \in \Sigma \, / \, u\geq 0\}, \quad \Sig=\{u\in H^1(S^4_+)/\, \, ||u||=1\}.
\end{eqnarray}
The Palais-Smale condition fails to be satisfied for $J$ on $\Sigma ^+$.
To characterize  the sequences failing the Palais-Smale condition, we
need to fix some notation.\\
For $a\in \overline{S_+^4}$ and $\l >0$, let
\begin{eqnarray}\label{e:o3}
\d _{a,\l }(x)=\frac {2\sqrt{2}\l}{\l ^2+1+(1-\l ^2) \cos d(a,x)},
\end{eqnarray}
where $d$ is the geodesic distance on $(\ov{S^4_+},g)$.
This function satisfies the following equation
$$
-\D \d _{a,\l } + 2\d _{a,\l }=\d _{a,\l }^3, \quad \mbox{ in }
S_+^4.
$$
Let $\var_{(a,\l)}$ be the function defined on $ S_+^4$ and satisfying
$$
\D \var _{(a,\l) } + 2\var _{(a,\l) }=\D \d _{a,\l } + 2\d _{a,\l }  \mbox{ in } S_+^4, \quad \frac{\partial\var_{(a,\l)}}{\partial\nu}=0  \mbox{ on } \partial S_+^4.
$$
Now, for $\e >0$ and $p\in \N^*$, let us define
\begin{align*}
V(p,\e )=& \{u\in \Sig /\exists a_1,...,a_p \in \overline{S_+^4}, \exists
\l _1,...,\l _p >0, \exists \a _1,...,\a _p>0
 \mbox{ s.t. } ||u-\sum_{i=1}^p\a _i\d _i||<\e, \\
 &  |\frac
{\a _i ^2K(a_i)}{\a _j ^2K(a_j)}-1|<\e ,
 \l _i>\e ^{-1}, \e _{ij}<\e \mbox{ and }\l _id_i<\e \mbox{ or }
\l _id_i >\e ^{-1}\},
\end{align*}
where $\d _i=\d _{a_i,\l _i}$, $d_i=d(a_i,\partial S_+^4)$ and
$\e _{ij}^{-1}=\l _i/\l _j +\l _j/\l _i + \l _i\l _j(1-\cos d(a_i,a_j))/2$.\\ \\
The failure of the Palais-Smale condition can be described, following the ideas introduced in \cite{BC}, \cite{L},  \cite{S} as follows:
\begin{pro}\label{p:22} 
Assume that $J$ has no critical point in $\Sigma ^+$ and let $(u_k)\in
\Sigma ^+$ be a sequence such that $J(u_k)$ is bounded and $\n J(u_k)\to 0$.
Then there exist an integer $p\in \N^*$, a sequence $\e _k>0$ ($\e _k\to 0$)
and an extracted subsequence of $u_k$, again denoted $(u_k)$, such that
$u_k\in V(p,\e _k )$.
\end{pro}
If a function $u$ belongs to $V(p,\e)$, we assume that, for the sake of simplicity, $\l_id_i < \e$ for $i\leq q$ and $\l_id_i > \e^{-1}$ for $i > q$. We consider the following minimization problem for $u\in V(p,\e)$ with $\e$ small
\begin{eqnarray}\label{e:51}
\min\{||u-\sum_{i=1}^q\a _i\d _{(a_i,\l_i)}-\sum_{i=q+1}^p\a _i\var_{(b_i,\l_i)} ||,\, \a _i>0,\, \l _i>0,\, a_i\in
\partial S_+^4 \mbox{ and } b_i\in S_+^4\}.
\end{eqnarray}
We then have the following proposition which defines a parametrization of the set $V(p,\e )$. It follows from the corresponding  statements in \cite{B2}, \cite{BC'}, \cite{R}.

\begin{pro}\label{p:23}

For any $p\in \N^*$, there is $\e _p>0$ such that if $\e <\e _p$ and $u\in V(p,\e )$, the minimization problem \eqref{e:51}
has a unique solution (up to permutation). In particular, we can write
$u\in V(p,\e )$ as follows
$$
u=\sum_{i=1}^q\bar{\a }_i\d _{(\bar{a}_i,\bar{\l }_i)}+\sum_{i=q+1}^p\bar{\a }_i\var _{(\bar{a}_i,\bar{\l }_i)}+ v,
$$
where $(\bar{\a }_1,...,\bar{\a }_p,\bar{a}_1,...,\bar{a}_p,\bar{\l }_1,...,
\bar{\l }_p)$ is the solution of \eqref{e:51} and $v\in H^1(S_+^n)$ such that
$$(V_0)\qquad
(v,\psi)=0 \mbox{ for } \psi\in \bigg\{\d_i,\frac{\partial\d_i}{\partial\l_i},\frac{\partial\d_i}{\partial a_i},\var_j,\frac{\partial\var_j}{\partial\l_j},\frac{\partial\var_j}{\partial a_j}/ \, \, i\leq q , \, j>q\bigg\}.
$$
\end{pro}
We also have the following proposition whose proof is similar, up to minor modifications to the corresponding statements in \cite{B1} (see also \cite{R})
\begin{pro}\label{p:24}
There exists a $C^1$ map which, to each $(\a_1,...,\a_ p,a_1,...,a_ p,\l_1,..., \l_ p)$ such that
$ \sum_{i=1}^p\a _i\d _i \in V(p,\e )$ with small $\e $,
associates $\ov{v}=\ov{v}_{(\a_i,a_i,\l_i )}$ satisfying
$$
J\left(\sum_{i=1}^q\a _i\d _i +\sum_{i=q+1}^p\a_i\var_i+\ov{v}\right)= \min\bigg\{
J\left( \sum_{i=1}^q\a _i\d _i+\sum_{i=q+1}^p\a_i\var_i +v\right) , \, v \mbox{ satisfies } (V_0)\bigg\}.
$$
Moreover, there exists $c>0 $ such that the following holds
$$
||\ov{v}||\leq c \left(\sum_{i\leq q}\frac{1}{\l
_i}+\sum_{i>q}\frac{|\n
K(a_i)|}{\l_i}+\sum_{i>q}\frac{1}{(\l_id_i)^2}+\sum_{k\ne r}\e
_{kr} (log (\e _{kr}^{-1}))^{1/2}\right).
$$
\end{pro}
Next we are going to recall a useful expansion of functional $J$ in $V(p,\e)$.
\begin{pro}\label{p:25} \cite{BEO2}
For $\e>0$ small enough and
$u=\sum_{i=1}^p\a_i\var_{(a_i,\l_i)}\in V(p,\e)$, we have the
following expansion
\begin{align*}
J(u)= & \frac{S_4^{1/2}\sum\a_i ^2}{\left(\sum\a_i ^4K(a_i)\right)^{1/2}}\left(1+\frac{\o_3}{8S_4}\left(\sum K(a_i)^{-1}\right)^{-1}\left(\sum\biggl(\frac{-\D K(a_i)}{3\l_i ^2K(a_i)^2}+\frac{4H(a_i,a_i)}{\l_i ^2K(a_i)}\biggr)\right.\right.\\
 &\left.\left. -\sum_{i\ne j}\frac{2}{\left(K(a_i)K(a_j)\right)^{1/2}}\left(\e_{ij}-\frac{2H(a_i,a_j)}{\l_i\l_j}\right)\right)+o\left(\sum \e_{kr}+\frac{1}{(\l_kd_k)^2}\right)\right),
\end{align*}
where $S_4 = 64\int_{\R^4}\frac{dx}{(1+|x|^2)^4}$.
\end{pro}

\section{Proof of Theorems}
\mbox{} Before  giving the proof of our theorems, we extract from
\cite{BEO2} the characterization of   the critical points at
infinity of our problem. We recall that the critical points at
infinity are  the orbits of the gradient flow of $J$ which remain
in $V(p,\e(s))$, where $\e(s)$ , a given function, tends to zero
when $s$ tends to $+\infty$ (see \cite{B1}).
 \begin{pro}\label{p:31}  \cite{BEO2}
Assume that for any $s$, $M(\tau_s)$ is nondegenerate. Thus, for $p\geq 1$, there exists a pseudogradient $W$ so that the
following holds:\\ There is a constant $c > 0$ independent of
$u=\sum_{i=1}^q\a_i\d_i+\sum_{j=q+1}^p\a_j\var_j\in V(p,\e)$ so
that
$$
(-\n J(u),W)\geq c\biggl(\sum_{k\ne r}\e_{kr}+ \sum_{i\leq
q}\frac{1}{\l_i}+\sum_{j=q+1}^p \frac{|\n K(a_j)|}{\l_j}+
\frac{1}{(\l_jd_j)^2}\biggr)\leqno{ (i)}
$$
$$
(-\n J(u+\ov{v}),W+\frac{\partial\ov{v}} {\partial
(\a_i,a_i,\l_i)}(W))\geq c\biggl(\sum_{k\ne r}\e_{kr}+ \sum_{i\leq
q}\frac{1}{\l_i}+ \sum_{j=q+1}^p \frac{|\n
K(a_j)|}{\l_j}+\frac{1}{(\l_jd_j)^2}\biggr)\leqno{ (ii)}
$$
(iii) $|W|$ is bounded. Furthermore, the only case where the
maximum of the $\l_i$'s is not bounded is when each point $a_i$ is
close to a critical point $y_{j_i}$ of $K$ with $j_i\ne j_k$ for
$i\ne k$ and $\rho(y_{i_1},...,y_{i_p})>0$, where
$\rho(y_{i_1},...,y_{i_p})$ denotes the least eigenvalue of
$M(y_{i_1},...,y_{i_p})$.
\end{pro}
\begin{cor}\label{c:32} \cite{BEO2}
Assume that for any $s$ $M(\tau_s)$ is nondegenerate, and assume further that $J$ has no critical point in $\Sig^+$. Then the only
critical points at infinity of $J$ correspond to
$$
\sum_{j=1}^pK(y_{i_j})^{-1/2}\var _{(y_{i_j},\infty)}, \quad\mbox{
with }\,  p\in \N^* \mbox{ and }\rho(y_{i_1},...,y_{i_p}) > 0.
$$
In addition, in the neighborhood of such a critical point at
infinity, we can find a change of variable
$$
(a_1,...,a_p,\l_1,...,\l_p)\to (\tilde a_1,...,\tilde
a_p,\tilde\l_1,...,\tilde\l_p):=(\tilde a,\tilde\l)
$$
such that
$$
J\biggl(\sum_{i=1}^p\a_i\var_i+\ov{v}\biggr)=
\psi(\a,\tilde a,\tilde\l):=\frac{8S_4^{1/2}\sum_{i=1}^p\a_i ^2}
{(\sum_{i=1}^p\a_i ^4K(a_i))^{1/2}}\biggl(1+(c-\eta)
\biggl(\sum_{i=1}^p\frac{1}{K(y_{j_i})}\biggr)^{-1}\, {}^t\L M(\tau_p)\L\biggr)
$$
where $\a=(\a_1,...,\a_p)$, $c$ is a positive constant, $\eta$ is
 a small positive constant, ${}^t\L=(\tilde\l_1,...,\tilde\l_p)$,
 $\tau_p=(y_{j_i},...,y_{j_p})$.\\
\end{cor}
Now we are ready to prove our theorems.\\
\begin{pfn}{\bf Theorem \ref{t:11} }
For $\eta > 0$ small enough, we introduce, following \cite{BEO2},
this neighborhood of $\Sig^+$ $$ V_\eta(\Sig^+)=\{u\in \Sig/\,
\,e^{2J(u)} J(u)^3|u^-|_{L^4}^2 < \eta \},
$$
where $u^-=\max(0,-u)$.\\
Recall that, from Proposition \ref{p:31} we have a vector field
$W$ defined in $V(p,\e)$ for $p\geq 1$. Outside $\cup_{p\geq
1}V(p,\e/2)$, we will use $-\n J$ and our global vector field $Z$
will be built using a convex combination of $W$ and $-\n J$.
$V_\eta(\Sig^+)$ is invariant under the flow line generated by $Z$
(see \cite{BCCH}). Arguing by contradiction, we assume that $J$
has no critical point in $V_\eta(\Sig^+)$. For any $y$ critical
point of $K$, set
$$
c_\infty (y) = \biggl (\frac{S_4}{K(y)} \biggr )^{1/2}.
$$
Since $y_0$ is the unique absolute maximum of $K$, we derive that
$$
c_\infty (y_0) < c_\infty (y), \qquad \forall y \neq y_0,
$$
where $y$ is any critical point of $K$.\\
 Let $u_0 \in \Sig^+$ such that
\begin{eqnarray}\label{e:31}
c_\infty (y_0) < J(u_0) < \inf _{y/ y\neq y_0, \n K(y) =0}c_\infty (y)
\end{eqnarray}
and let $\eta (s,u_0)$ be the one parameter group generated by
$Z$. It is known that $|\n J|$ is lower bounded outside $V(p,\e
/2)$, for any $p\in \N^*$ and for $\e$ small enough, by a fixed
constant which depends only on $\e$. Thus the flow line $\eta
(s,u_0)$ cannot remain outside of the set  $V(p,\e /2)$.
Furthermore, if the flow line travels from  $V(p,\e /2)$ to the
boundary of  $V(p,\e )$, $J(\eta (s,u_0))$ will decrease by a
fixed constant which depends on $\e$. Then, this travel cannot be
repeated in an infinite time. Thus there exist $p_0$ and $s_0$
such that the flow line enters into  $V(p_0,\e /2)$ and it does
not exit from $V(p_0,\e)$. Since $u_0$ satisfies \eqref{e:31}, we
derive that $p_0=1$, thus, for $s\geq s_0$,
$$
\eta (s,u_0)=\a _1 \var_{ (x_1(s),\l _1(s))} + v(s)
$$
Using again \eqref{e:31}, we deduce that $x_1(s)$ is outside
$\mathcal{V}(y,\tau )$ for any $y \in \mathcal{F}^+ \diagdown
\{y_0\}$, where $\mathcal{V}(y,\tau )$ is a neighborhood of $y$
and where $\tau$ is a small positive real. Now, by  assumptions
of Theorem \ref{t:11} and by the construction of a pseudogradient
$Z$, we derive that $\l _1(s)$ remains bounded along the flow
lines of $Z$. Thus we obtain
$$
|\n J(\eta (s,u_0)).Z(\eta (s,u_0))| \geq c > 0 \quad \forall s\ge 0,
$$
where $c$ depends only on $u_0$.\\
 Then when $s$ goes to $+\infty$, $J(\eta (s,u_0))$ goes to
$-\infty$ and this yields a contradiction. Thus there exists a
critical point of $J$ in $V_\eta(\Sig^+)$. Arguing as in
\cite{BCCH}, we prove that such a critical point is positive and
hence our result follows.
\end{pfn}\\
Now before giving the proof of Theorem \ref{t:12}, we state the
following lemma. Its proof is very similar to the proof of
Corollary B.3 of \cite{BC'}(see also \cite{B2}), so we will omit
it.
\begin{lem}\label{l:33}
Let $a_1$, $a_2 \in S^4_+$, $\a_1$, $\a_2 > 0$ and $\l$ large
enough. For  $u=\a_1\var_{ (a_1,\l )} +\a_2\var_{ (a_2,\l )}$, we
have
$$
J\biggl(\frac{u}{||u||}\biggr) \leq
S_4^{1/2}\biggl(\frac{1}{K(a_1)} + \frac{1}{K(a_2)}
\biggr)^{1/2}(1+o(1)).
$$
\end{lem}
\begin{pfn}{\bf Theorem \ref{t:12} }
Again, we argue by contradiction. We assume that $J$ has no
critical point in $V_\eta (\Sig^+)$. Let
$$
c_\infty (y_0,y_{i_1}) = S_4^{1/2}\biggl(\frac{1}{K(y_0)} +
\frac{1}{K(y_{i_1})} \biggr)^{1/2}
$$
We observe that under the assumption $(A_1)$ of Theorem
\ref{t:12}, $(y_0,y_{i_1})$ is not a critical point at infinity of
$J$. Using Corollary \ref{c:32} and the assumption $(A_2)$ of
Theorem \ref{t:12}, it follows that the only critical points at
infinity of $J$ under the level $c_1=c_\infty (y_0,y_{i_1}) + \e$,
for $\e$ small enough, are $\var_{ (y_0,\infty )}$ and   $\var_{
(y_{i_1},\infty )}$. The unstable manifolds at infinity of such
critical points at infinity, $W_u(y_0)_\infty$,
$W_u(y_{i_1})_\infty$ can be described, using Corollary
\ref{c:32}, as the product of $W_s(y_0)$,
 $W_s(y_{i_1})$ (for a pseudogradient of $K$) by $[A, +\infty [$
 domain of the variable $\l$, for some positive
number $A$ large enough.\\
Since $J$ has no critical point, it follows that $ J_{c_1}=\{u\in
\Sig ^+ / J(u) \leq c_1 \}$ retracts by deformation on $X_\infty =
 W_u(y_0)_\infty \cup W_u(y_{i_1})_\infty$ (see Sections 7 and 8 of
\cite{BR}) which can be parametrized  by $X \times
[A, +\infty[$, where $X=\overline{W_s(y_{i_1})}$.\\
Under (1.3) and (1.4) (see the first section), $X= W_s(y_0) \cup
W_s(y_{i_1})$. Thus $X$ is a manifold in dimension $k_{i_1}$ without
boundary.\\
We claim now that $X_\infty$ is contractible in $J_{c_1}$. Indeed,
let $h: [0,1] \times X \times \left[A,\right. +\infty\left[
\right.\longmapsto \Sigma ^+$ defined by
$$
 (t,x,\l) \longmapsto \frac{t\var _{(y_0,\l )} +
 (1-t)\var _{(x,\l )}}{||t\var _{(y_0,\l )} + (1-t)\var _{(x,\l )}||}
$$
$h$ is continuous and satisfies
$$
h(0,x,\l )=\frac{\var _{(x,\l )}}{||\var _{(x,\l )}||} \quad
\mbox{and}\quad  h(1,x,\l )=\frac{\var _{(y_0,\l )}}{||\var
_{(y_0,\l )}||}.
$$
In addition, since $K(x)\geq K(y_{i_1})$ for any $x\in X$, it
follows from Lemma 3.3 that $J(h(t,x,\l))<c_1$, for each
$(t,x,\l)\in [0,1] \times X \times \left[A,\right. +\infty\left[
\right.$. Thus the contraction $h$ is performed under the level
$c_1$. We derive that $X_\infty$ is contractible in $J_{c_1}$,
which retracts by deformation on $X_\infty$, therefore $X_\infty$
is contractible leading to the contractibility of $X$, which is a
contradiction, since $X$ is a manifold in dimension $k_{i_1}$ without
boundary. Hence there exists a critical point of $J$ in $V_\eta
(\Sigma ^+ )$. Arguing as in \cite{BCCH}, we prove that such a
critical point is positive. Now we are going
to show that such a critical point has a Morse index equal to $k_{i_1}$ or $k_{i_1}+1$.\\
Using a dimension argument and since $h([0,1], X_\infty )$ is a
manifold in dimension $k_{i_1}+1$, we derive that the Morse index of
such a critical point is $\leq k_{i_1}+1$.\\
Now, arguing by contradiction, we assume that the Morse index is $\leq k_{i_1}-1$.
Perturbing, if necessary $J$, we may assume that all the critical
points of $J$ are nondegenerate and have their Morse index $\leq
k_{i_1}-1$.\\
 Such critical points do not change the homological group in
dimension $k_{i_1}$ of level sets of $J$.\\
Now let  $c_\infty (y_{i_1})= S_4^{1/2}K(y_{i_1})^{-1/2}$ and let
$\e$ be a small positive real. Since $X_\infty$ defines a
homological class in dimension $k_{i_1}$ which is trivial in $J_{c_1}$,
but not trivial in  $J_{c_\infty (y_{i_1})+ \e}$, our result
follows.
\end{pfn}\\
\begin{pfn}{\bf Theorem \ref{t:13}}
We notice that the assumption $(H_3)$ implies that $(y_0,y_{i_1})$
is a critical point at infinity of $J$. Now, arguing by
contradiction, we assume that $(1)$ has no solution. We claim that
$f_\l (C_{y_o}(X))$ retracts by deformation on $X\cup
W_u(y_0,y_{i_1})_\infty$. Indeed, let
$$
u=\a\var _{(y_0,\l )} + (1-\a )\var _{(x,\l )} \in f_\l (C_{y_0}(X)),
$$
the action of the flow of the pseudogradient $Z$ defined in
the proof of Theorem \ref{t:11} is essentially on $\a$ (see [5] and [11]).\\
- If $\a < 1/2$, the flow of $Z$ brings $\a$ to zero and
thus $u$ goes in this case to $\overline{W_u(y_0)}_\infty \equiv \{y_0\}$.\\
- If  $\a > 1/2$, the flow of $Z$ brings $\a$ to $1$ and thus $u$ goes,
in this case, to $\overline{W_u(y_{i_1})}_\infty \equiv X_\infty$.\\
- If $\a = 1/2$, since only $x$ can move then $y_0$ remains one of
the points of concentration of $u$ and $x$ goes to $W_s(y_i)$, where
$y_i=y_{i_1}$ or $y_i=y_0$ and two cases may occur :\\
- In the first case, that is, $y_i=y_{i_1}$, $u$ goes to $W_u(y_0, y_{i_1})_\infty$ .\\
- In the second case, that is, $y_i=y_0$, there exists $s_0 \geq
0$ such that $x(s_0)$ is close to $y_0$.  Thus, using Lemma
\ref{l:33}, we have the following inequality
$$
J(u(s_0)) \leq c_\infty (y_0,y_0) + \gamma := c_2,
$$
where $c_\infty (y_0,y_0)= S_4^{1/2} (2/K(y_0))^{1/2}$ and
where $\gamma$ is a positive constant small enough.\\
Now, using  assumption $(H_4)$, it follows from Corollary
\ref{c:32} that $J_{c_2}$ retracts by deformation on
$W_u(y_0)_\infty \equiv \{y_0\}$ and thus $u$ goes
to $W_u(y_0)_\infty$. Therefore $f_\l (C_{y_0}(X))$
retracts by deformation on $X_\infty \cup W_u(y_0,y_{i_1})_\infty$.\\
Now, since $\mu (y_{i_1}) =0$, it follows that this strong retract
does not intersect $W_u(y_0,y_{i_1})_\infty$ and thus it is
contained in $X_\infty$. Therefore $X_\infty$ is contractible,
leading to the contractibility of $X$, which is a contradiction,
since $X$ is a manifold of dimension $k_{i_1}$ without boundary. Hence
$(1)$ admits a solution. Now, using the same arguments as those
used in the proof of Theorem \ref{t:12}, we easily derive that the
Morse index of the solution provided above is equal to $k_{i_1}$ or
$k_{i_1}+1$. Thus our result follows.
\end{pfn}\\
\begin{pfn}{\bf Theorem \ref{t:14}}
Assume that $(1)$ has no solution. Using the same arguments as
those in the proof of Theorem \ref{t:13}, we deduce that  $f_\l
(C_{y_0}(X))$ retracts by deformation on $$ X_\infty \cup
\left(\cup _{y\in B_k}W_u(y_0, y)_\infty \right)\cup D,
$$
where $D\subset \sigma$ is a stratified set and where $\sigma =
\cup _{y\in X\diagdown B_k}W_u(y_0, y)_\infty$
is a manifold in dimension at most $k$.\\
Since $\mu  (y)=0$ for each $y \in B_k$,  $f_\l (C_{y_0}(X))$
retracts by deformation on $X_\infty \cup D$, and therefore
$H_*(X_\infty\cup D)=0$, for all $*\in \N^*$, since
$f_\l(C_{y_0}(X))$ is a contractible set. Using the exact homology
sequence of $(X_\infty\cup D, X_\infty)$, we obtain
$$
...\to H_{k+1}(X_\infty\cup D) \to ^{\pi}  H_{k+1}(X_\infty\cup D,
X_\infty) \to ^{\partial}  H_k(X_\infty) \to ^{i}
H_k(X_\infty\cup D) \to ...
$$
Since $H_*(X_\infty\cup D) = 0$, for all $*\in \N^*$, then
$H_k(X_\infty)=H_{k+1}(X_\infty\cup D, X_\infty)$.\\
In addition, $(X_\infty\cup D, X_\infty)$ is a stratified set of
dimension at most $k$, then $H_{k+1}(X_\infty\cup D,X_\infty) =
0$, and therefore $H_k(X_\infty)=0$. This implies that $H_k(X)=0$
(recall that $X_\infty\equiv X\times [A,\infty)$). This yields a
contradiction since $X$ is a manifold in dimension $k$ without
boundary. Then, arguing as in the end of the proof of Theorem
\ref{t:13}, our theorem follows.
\end{pfn}\\
\begin{pfn}{\bf Theorem \ref{t:15} }
We argue by contradiction. Assume that $(1)$ has no solution.
 Let
$c_1=\frac{3}{2}S_4^{1/2}$.  Using the expansion of $J$, see
Proposition \ref{p:25}, we derive that there exists a constant
$c_0$ independent of $K$ such that if $||K-1||_{L^\infty (S^4_+)}
\leq c_0$, then the following holds
$$
J(u)< c_1, \quad \forall u\in V(1,\e ) \quad \mbox{and} \quad J(u)
> c_1 \quad \forall u \in V(p,\e ) \,\,\mbox{with }\,\, p\geq 2,
$$
where $\e$ is a small positive real.\\
Therefore, it follows from Corollary \ref{c:32} that the critical
points at infinity of $J$ under the level $c_1$ are in one to one
correspondence with the critical points $y$ of $K$ such that $y\in
\mathcal{F}^+$. Since $J$ has no critical points in $\Sig ^+$, it
follows that $J_{c_1}$ retracts by deformation on $X_\infty = \cup
_{y_i\in\mathcal{F}^+}W_u(y_i)_\infty$ (see sections 7 and 8 of
\cite{BR}) which can be parametrized by
 $X \times [A, +\infty[$ (we recall that $X=\cup_{y\in \mathcal{F}^+} \ov{W_s(y)}$).\\
Now we claim that $X_\infty$ is contractible in
$J_{c_1}$. Indeed, since $S^4_+$ is a contractible set, we deduce
that there exists a contraction $h: [0,1] \times X \to S^4_+$, $h$
continuous and satisfies for any $a\in X$, $h(0,a)=a$ and
$h(1,a)=a_0$ a point of $S^4_+$. Such a contraction gives rise to
the following contraction $\tilde{h}:  [0,1] \times X_\infty \to
\Sig ^+$ defined by
$$
[0,1] \times X \times \left[A,\right.+\infty\left[ \right.\ni
(t,a,\l  )\longmapsto \var _{(h(t,a),\l )} + \bar{v} \in \Sigma
^+.
$$
For $t=0$, $\var _{(h(0,a_),\l )}+\bar{v} = \var _{(a, \l ) }
+\bar{v} \in X_\infty$. $\tilde{h}$ is continuous and
$\tilde{h}(1,a,\l )= \var _{(a_0,\l )} +\bar{v}$, hence our
claim follows.\\
Now, using Proposition \ref{p:25}, we deduce that
$$
J(\var _{(h(t,a), \l) } + \bar{v}) \sim
S_4^{\frac{1}{2}}(K(h(t,a)))^{\frac{-1}{2}}\left(1+O(A
    ^{-2})\right).
$$
\\
Choosing $c_0$ small enough and $A$ large enough, we then have
$J(\var _{(h(t,a),\l)}) < c_1$. Therefore such a contraction is
performed under the level $c_1$, so $X_\infty$ is contractible in
$J_{c_1 }$, which retracts by deformation on $X_\infty$, therefore
$X_\infty$ is contractible leading to the contractibility of $X$,
which is in contradiction with assumption $(H_{8})$. We derive
that (1) has a solution. Lastly, the same arguments as those in the proof of Theorem \ref{t:12} easily give the desired estimate for the Morse index of the solution provided above  and therefore our result follows.
\end{pfn}\\
\begin{pfn}{\bf Theorem \ref{t:16} }
Arguing by contradiction, we assume that $(1)$ has no solution. Notice that under the assumption of our theorem,  $J$ has no critical point at
infinity having more or equal to two masses in its description and
therefore the only critical points at infinity of $J$ are
$\var(y,\infty)$ such that $y\in \mathcal{F}^+$. Since $J$ has no
critical point, it follows that $\Sig^+$ retracts by deformation
on $X_\infty=\cup_{y\in \mathcal{F}^+}W_u(y)_\infty$, see
\cite{BR}. Thus we conclude that $X_\infty$ is contractible since
$\Sig^+$ is a contractible set. So $X$ is contractible and this is
a contradiction with our assumption. Hence $(1)$ has a solution and as above, we derive that the Morse index of such a solution is $\geq$ $m$, therefore  our result follows.
\end{pfn}\\

\end{document}